\documentclass[11pt,twoside]{amsart}

\usepackage{amsmath}
\usepackage{amsfonts}
\usepackage{amssymb}
\usepackage{amsthm}
\usepackage{color}
\usepackage{graphicx}
\usepackage[all]{xy}
\usepackage{fullpage}
\usepackage{url}

\theoremstyle{definition}
\newtheorem{defn}{Definition}[section]

\theoremstyle{plain}
\newtheorem{thm}[defn]{Theorem}

\newtheorem{conj}[defn]{Conjecture}

\def\P{\ensuremath{\mathbb{P}}}

\def\R{\ensuremath{\mathbb{R}}}

\def\FF{\ensuremath{\mathcal F}}

\def\OO{\ensuremath{\mathcal O}}

\def\TT{\ensuremath{\mathcal T}}

\def\ch{\mathop{\mathrm{ch}}\nolimits}

\def\Coh{\mathop{\mathrm{Coh}}\nolimits}

\def\into{\ensuremath{\hookrightarrow}}
\def\onto{\ensuremath{\twoheadrightarrow}}

\begin{document}

\title[Counterexample to the Generalized Bogomolov-Gieseker Inequality]{Counterexample to the Generalized Bogomolov-Gieseker Inequality for Threefolds}

\author{Benjamin Schmidt}
\address{Department of Mathematics, The Ohio State University, 231 W 18th Avenue, Columbus, OH 43210-1174, USA}
\email{schmidt.707@osu.edu}
\urladdr{https://people.math.osu.edu/schmidt.707/}



\begin{abstract}
We give a counterexample to the generalized Bogomolov-Gieseker inequality for threefolds conjectured by Bayer, Macr\`i and Toda using the blow up of a point over three dimensional projective space.
\end{abstract}

\maketitle

\section{Introduction}

Tom Bridgeland introduced the notion of a stability condition on a triangulated category in \cite{Bri07}. He was inspired by the study of Dirichlet branes in string theory by Douglas in \cite{Dou02}. Instead of defining stability such as slope stability or Gieseker stability in the category of coherent sheaves, one uses other abelian categories inside the bounded derived category of coherent sheaves. One of the main problems is the construction of stability conditions on smooth projective threefolds.

In \cite{BMT14} Bayer, Macr\`i and Toda proposed a conjectural construction for threefolds based on a generalized Bogomolov-Gieseker inequality involving the third Chern character. It was proven for $\P^3$ in \cite{MacE14}, for the smooth quadric hypersurface $Q$ in $\P^4$ in \cite{Sch14} and for abelian threefolds independently in \cite{MP13b} and \cite{BMS14}. Moreover, the cases of $\P^3$ and $Q$ were recently generalized to all Fano threefolds of Picard rank $1$ in \cite{Li15}. 

Let $X$ be the blow up of $\P^3$ in one point. The article \cite{Li15} contained the conjecture for $X$ as a further example. Unfortunately that example contained a mistake, which lead us directly to the counterexample. In order to do so, we choose the polarization to be $H = -1/2 K_X$ and the pullback of $\OO_{\P^3}(1)$ is the object that does not satisfy the conjecture.

\subsection*{Acknowledgements}

The author would like to thank Arend Bayer, Chunyi Li, Emanuele Macr\`i, Yukinobu Toda, Bingyu Xia and Xiaolei Zhao for commenting on preliminary versions of this article. The author is partially supported by NSF grant DMS-1523496 (PI Emanuele Macr\`i) and a Presidential Fellowship of the Ohio State University.

\section{Preliminaries}

We start by recalling the notion of tilt stability due to \cite{BMT14}. Let $X$ be a smooth projective threefold over the complex numbers and let $H$ be an ample divisor on $X$. For $\alpha \in \R_{> 0}$, $\beta \in \R$ and $E \in D^b(X)$ we define
\[
\nu_{\alpha, \beta}(E) = \frac{H \cdot \ch_2^{\beta}(E) - \frac{\alpha^2}{2} H^3 \cdot \ch_0^{\beta}(E)}{H^2 \cdot \ch_1^{\beta}(E)},
\]
where $\ch^{\beta}(E) = e^{-\beta H} \cdot \ch(E)$. A torsion pair is defined by
\begin{align*}
\TT_{\beta} &= \{E \in \Coh(X) : \forall E \onto G: H^2 \cdot \ch_1(G) > \beta H^3 \cdot \ch_0(G) \text{ or } H^2 \cdot \ch_1(G) = \ch_0(G) = 0 \}, \\
\FF_{\beta} &=  \{E \in \Coh(X) : \forall F \subset E: H^2 \cdot \ch_1(G) \leq \beta H^3 \cdot \ch_0(G) \text{ and } \ch_0(G) \neq 0 \}.
\end{align*}
A new heart of a bounded t-structure is then defined as the extension closure $\Coh^{\beta}(X) = \langle \FF_{\beta}[1], \TT_{\beta} \rangle$.

\begin{defn}
An object $E \in \Coh^{\beta}(X)$ is $\nu_{\alpha, \beta}$-\textit{semistable} if for all subobjects $F \into E$ the inequality $\nu_{\alpha, \beta}(F) \leq \nu_{\alpha, \beta}(E)$ holds.
\end{defn}

Bayer, Macr\`i and Toda conjectured an inequality for tilt semistable objects involving all three Chern characters in \cite{BMT14}. Over time their conjecture has evolved into the following form. We define 
\begin{align*}
\Delta(E) &= (H^2 \cdot \ch_1(E))^2 - 2 (H^3 \cdot \ch_0(E)) (H \cdot \ch_2(E)), \\
Q_{\alpha, \beta}(E) &= \alpha^2 \Delta(E) + 4 (H \cdot \ch_2^{\beta}(E))^2 - 6 (H^2 \cdot \ch_1^{\beta}(E)) \ch_3^{\beta}(E).
\end{align*}

\begin{conj}[\cite{BMT14, BMS14}]
Any $\nu_{\alpha, \beta}$-semistable object $E \in \Coh^{\beta}(X)$ satisfies $Q_{\alpha, \beta}(E) \geq 0$.
\end{conj}

In order to prove the counterexample in the next section, we need to have a better understanding of the structure of walls in tilt stability. A \textit{numerical wall} is a non trivial zero set of a non trivial equation of the form $\nu_{\alpha, \beta}(v) = \nu_{\alpha, \beta}(w)$ for some $v,w \in K_0(X)$. Numerical walls in tilt stability satisfy Bertram's Nested Wall Theorem. For surfaces it was proven in \cite{MacA14}. It still holds in the threefolds case as for example shown in \cite{Sch15}.

\begin{thm}[Structure Theorem for Walls in Tilt Stability]
Let $v \in K_0(X)$ be fixed. All numerical walls in the following statements are with respect to $v$.
\begin{enumerate}
  \item Numerical walls in tilt stability are of the form 
  \[x\alpha^2 + x\beta^2 + y\beta + z = 0,\]
  for $x,y,z \in \R$, where $x = 0$ if and only if both classes defining the wall have the same classical slope. In particular, they are either semicircles with center on the $\beta$-axis or vertical rays.
  \item Numerical walls only intersect if they are identical.
  \item If a numerical wall has a single point at which it is an actual wall, then all of it is an actual wall.
  \item The equation $Q_{\alpha, \beta}(v) = 0$ is a numerical wall for $v$.
\end{enumerate}
\end{thm}

\section{Counterexample}

Let $f: X \to \P^3$ be the blow up of $\P^3$ in a point $P$. The Picard group of $X$ is well known to be a free abelian group with two generators $\OO(L) = f^*\OO_{\P^3}(1)$ and $\OO(E)$, where $E = f^{-1}(P)$. The variety $X$ is Fano with  canonical divisor given by $-4L + 2E$. In particular, $H = 2L - E$ is an ample divisor. We also have the intersection products $L^3 = E^3 = 1$, $L \cdot E = 0$ and $H^3 = 7$. The goal of this section is to prove the following counterexample to the conjectural inequality.

\begin{thm}
There exists $\alpha \in \R_{> 0}$ and $\beta \in \R$ such that the line bundle $\OO_X(L)$ is $\nu_{\alpha, \beta}$-stable, but $Q_{\alpha, \beta}(\OO_X(L)) < 0$.
\begin{proof}
Since $\OO(L)$ is a line bundle it is a slope semistable sheaf. In particular, either $\OO(L)$ or $\OO(L)[1]$ is $\nu_{\alpha, \beta}$-stable for all $\alpha \gg 0$.

A straightforward computation shows $H^3 \cdot \ch_0(\OO(L)) = 7$, $H^2 \cdot \ch_1(\OO(L)) = 4$, $H \cdot \ch_2(\OO(L)) = 1$ and $\ch_3(\OO(L)) = 1/6$. In particular, this means $H^2 \cdot \ch^{1/2}_1(\OO(L)) = 1/2$. If $F \into \OO(L)$ destabilizes along the line $\beta = 1/2$, we must have $\ch_1^{1/2}(F) \in \{0, 1/2\}$. That means either $F$ or the quotient has slope infinity independently of $\alpha$, a contradiction.

A completely numerical computation shows that $Q_{\alpha, \beta}(\OO(L)) \geq 0$ is equivalent to the inequality 
\[
\alpha^2 + \left(\beta - \frac{1}{4}\right)^2 \geq \frac{1}{16}.
\]
We are done if we can prove that there is no wall with equality in this inequality for $\OO(L)$. Assume there is a destabilizing sequence $0 \to F \to \OO(L) \to G \to 0$ giving exactly this wall. Taking the long exact sequence in cohomology, we get $H^3 \cdot \ch_0(F) \geq 7$. By definition of $\Coh^{\beta}(X)$ we have the inequalities $H^2 \cdot\ch_1^{\beta}(\OO(L)) \geq H^2 \cdot \ch_1^{\beta}(F) \geq 0$ for all $\beta \in [0, 1/2]$. This can be rewritten as
\[
4 + \beta (H^3 \cdot \ch_0(F) - 7) \geq H^2 \cdot \ch_1(F) \geq \beta H^3 \cdot \ch_0(F).
\]
Notice that the middle term is independent of $\beta$ and we can vary $\beta$ independently on the left and right. Therefore, we get $4 \geq H^2 \cdot \ch_1(F) \geq H^3 \cdot \ch_0(F)/2 $. This means $H^2 \cdot \ch_1(F) = 4$ and $H^3 \cdot \ch_0(F) = 7$. This does not give the correct wall.
\end{proof}
\end{thm}

\bibliography{all}

\def\cprime{$'$} \def\cprime{$'$}
\begin{thebibliography}{Mac14b}

\bibitem[BMS14]{BMS14}
A.~Bayer, E.~Macr{\`{\i}}, and P.~Stellari.
\newblock The space of stability conditions on abelian threefolds, and on some
  {C}alabi-{Y}au threefolds, 2014.
\newblock arXiv:1410.1585v1.

\bibitem[BMT14]{BMT14}
A.~Bayer, E.~Macr{\`{\i}}, and Y.~Toda.
\newblock Bridgeland stability conditions on threefolds {I}:
  {B}ogomolov-{G}ieseker type inequalities.
\newblock {\em J. Algebraic Geom.}, 23(1):117--163, 2014.

\bibitem[Bri07]{Bri07}
T.~Bridgeland.
\newblock Stability conditions on triangulated categories.
\newblock {\em Ann. of Math. (2)}, 166(2):317--345, 2007.

\bibitem[Dou02]{Dou02}
M.~R. Douglas.
\newblock Dirichlet branes, homological mirror symmetry, and stability.
\newblock In {\em Proceedings of the {I}nternational {C}ongress of
  {M}athematicians, {V}ol. {III} ({B}eijing, 2002)}, pages 395--408. Higher Ed.
  Press, Beijing, 2002.

\bibitem[Li15]{Li15}
C.~Li.
\newblock Stability conditions on fano threefolds of picard number one, 2015.
\newblock arXiv:1510.04089v2.

\bibitem[Mac14a]{MacA14}
A.~Maciocia.
\newblock Computing the walls associated to {B}ridgeland stability conditions
  on projective surfaces.
\newblock {\em Asian J. Math.}, 18(2):263--279, 2014.

\bibitem[Mac14b]{MacE14}
E.~Macr{\`{\i}}.
\newblock A generalized {B}ogomolov-{G}ieseker inequality for the
  three-dimensional projective space.
\newblock {\em Algebra Number Theory}, 8(1):173--190, 2014.

\bibitem[MP16]{MP13b}
A.~Maciocia and D.~Piyaratne.
\newblock Fourier--{M}ukai transforms and {B}ridgeland stability conditions on
  abelian threefolds {II}.
\newblock {\em Internat. J. Math.}, 27(1):1650007, 27, 2016.

\bibitem[Sch14]{Sch14}
B.~Schmidt.
\newblock A generalized {B}ogomolov-{G}ieseker inequality for the smooth
  quadric threefold.
\newblock {\em Bull. Lond. Math. Soc.}, 46(5):915--923, 2014.

\bibitem[Sch15]{Sch15}
B.~Schmidt.
\newblock Bridgeland stability on threefolds - some wall crossings, 2015.
\newblock arXiv:1509.04608v1.

\end{thebibliography}
\bibliographystyle{alphaspecial}

\end{document}